\documentstyle[12pt]{article}
\begin{document}
\title {Refinements of Reverse Triangle Inequalities in Inner Product Spaces}
\author{{\bf Arsalan Hojjat Ansari and Mohammad Sal Moslehian} \\ Dept. of Math., Ferdowsi Univ.\\ P. O. Box 1159, Mashhad 91775\\ Iran\\ E-mail: msalm@math.um.ac.ir\\URL: http://www.um.ac.ir/$\sim$moslehian/}
\date{}
\maketitle

\begin{abstract}
Refining some results of S. S. Dragomir, several new reverses of the triangle inequality in inner product spaces are obtained.
\end{abstract}
\newpage

\section{Introduction.}
It is interesting to know under which conditions the triangle inequality went the other way in a normed space $X$; in other words, we would like to know if there is a constant $c$ with the property that $c\sum_{k=1}^n\|x_k\|\leq\|\sum_{k=1}^n x_k\|$ for any finite set $x_1,\cdots,x_n\in X$. M. Nakai and T. Tada $\cite{N-T}$ proved that the normed spaces with this property are precisely those of finite dimensional.

The first authors investigating reverse of the triangle inequality in inner product spaces were J. B. Diaz and F. T. Metcalf $\cite{D-M}$ by establishing the following result as an extension of an inequality given by M. Petrovich $\cite{PET}$ for complex numbers:

{\bf Diaz-Metcalf Theorem.} Let $a$ be a unit vector in the inner product space $(H;\langle .,.\rangle)$. Suppose the vectors $x_k \in H, k\in\{1,\cdots,n\}$ satisfy
\begin{eqnarray*}
0\leq r\leq\frac{Re\langle x_k,a\rangle}{\|x_k\|}, k\in\{1,\cdots,n\}
\end{eqnarray*}
Then 
\begin{eqnarray*}
r\sum_{k=1}^n\|x_k\|\leq\|\sum_{k=1}^n x_k\|.
\end{eqnarray*}
where equality holds if and only if
\begin{eqnarray*}
\sum_{k=1}^n x_k =r\sum_{k=1}^n \|x_k\|a.
\end{eqnarray*}

Inequalities related to the triangle inequality are of special interest; cf. Chapter XVII of $\cite{M-P-F}$ and may be applied to get nice inequalities in complex numbers or to study vector-valued integral inequalities $\cite{DRA1}, \cite{DRA2}$.

Using several ideas and notation of $\cite{DRA1}, \cite{DRA2}$ we modify or refine results of S. S. Dragomir and get some new reverses of triangle inequality.

We use repeatedly the Cauchy-Schwarz inequality without mentioning it. The reader is refered to $\cite{RAS}, \cite{DRA3}$ for the terminology on inner product spaces.

\section{Main Results.}

The following theorem is an strengthen of theorem 1 of $\cite{DRA2}$ in which the real numbers $r_1, r_2$ are not neccesarily nonnegative. The proof seems to be different as well.

{\bf Theorem 1.} Let $a$ be a unit vector in the complex inner product space $(H;\langle .,.\rangle)$. Suppose the vectors $x_k \in H, k\in\{1,\cdots,n\}$ satisfy
\begin{equation}
0\leq r_1 ^2 \|x_k\|\leq Re\langle x_k,r_1 a\rangle, 0\leq r_2 ^2 \|x_k\| \leq Im\langle x_k,r_2 a \rangle
\end{equation}
for some $r_1,r_2\in[-1,1].$ Then we have the inequality
\begin{equation}
(r_1 ^2 +r_2 ^2)^\frac{1}{2}\sum_{k=1}^n\|x_k\|\leq\|\sum_{k=1}^n x_k\|.
\end{equation}
The equality holds in (2) if and only if
\begin{equation}
\sum_{k=1}^n x_k =(r_1 +ir_2)\sum_{k=1}^n \|x_k\|a.
\end{equation}

{\bf Proof.} If $r_1 ^2 +r_2 ^2 =0$,theorem is trivial. Assume that $r_1 ^2 +r_2 ^2 \neq 0$. Summing inequalities (1) over $k$ from $1$ to $m$, we have
\begin{eqnarray*}
(r_1 ^2 +r_2 ^2)\sum_{k=1}^n \|x_k\|\leq Re\langle \sum_{k=1}^n x_k ,r_1 a\rangle+Im\langle \sum_{k=1}^n x_k ,r_2a\rangle\\
=Re\langle \sum_{k=1}^n x_k ,(r_1+ir_2)a\rangle\\
\leq |\langle \sum_{k=1}^n x_k ,(r_1+ir_2)a\rangle|\\
\leq\|\sum_{k=1}^n x_k \| \|(r_1 +ir_2)a\|\\
=(r_1^2 +r_2^2)^\frac{1}{2}\|\sum_{k=1}^n x_k \|.
\end{eqnarray*}
Hence $2$ holds.

If (3) holds, then 
\begin{eqnarray*}
\| \sum_{k=1}^n x_k\|=\|(r_1 +ir_2)\sum_{k=1}^n \|x_k\|a\|=(r_1^2 +r_2^2)^\frac{1}{2}\sum_{k=1}^n \|x_k\|.
\end{eqnarray*}

Conversely, if the equality holds in (2), we have \\
\begin{eqnarray*}
(r_1^2 +r_2^2)^\frac{1}{2}\|\displaystyle{\sum_{k=1}^n x_k}\|&=&(r_1^2+r_2^2)\displaystyle{\sum_{k=1}^n \|x_k\|}\leq Re\langle\displaystyle{\sum_{k=1}^n x_k},(r_1+ir_2)a\rangle\\
&\leq&|\langle\displaystyle{\sum_{k=1}^n x_k},(r_1 +ir_2)a\rangle|\leq(r_1^2 +r_2^2)^\frac{1}{2}\|\displaystyle{\sum_{k=1}^n x_k\|}.
\end{eqnarray*}
From this we deduce 
$$|\langle\sum_{k=1}^n x_k ,(r_1 +ir_2)a\rangle|=\|\sum_{k=1}^n x_k\|\|(r_1 +ir_2)a\|.$$
Consequently there exists $\eta\geq 0$ such that $\displaystyle{\sum_{k=1}^n} x_k =\eta (r_1 +ir_2)a$

From this we have 
$$(r_1^2 +r_2^2)^\frac{1}{2}\eta=\|\eta(r_1 +ir_2)a\|=\|\sum_{k=1}^n x_k\|=(r_1^2 +r_2^2)^\frac{1}{2}\sum_{k=1}^n \| x_k\|.$$
Hence $\eta=\displaystyle{\sum_{k=1}^n \|x_k\|}.\Box$

The next theorem is a refinement of Corollary 1 of $\cite{DRA2}$ since, in the notation of the theorem, $\sqrt{2-p_1^2-p_2^2}\leq \sqrt{\alpha_1^2+\alpha_2^2}$.

{\bf Theorem 2.} Let $a$ be a unit vector in the complex inner product space $(H;\langle .,.\rangle)$. Suppose the vectors $x_k \in H-\{0\}, k\in\{1,\cdots,n\}$, 
such that
\begin{equation}
\|x_k-a\|\leq p_1, \|x_k-ia\|\leq p_2, p_1, p_2\in(0,\sqrt{\alpha^2+1})
\end{equation}
where $\alpha=\displaystyle{\min_{1\leq k\leq n}}\|x_k\|$. Let
\begin{eqnarray*}
\alpha_1=\min\{\frac{\|x_k\|^2-p_1^2+1}{2\|x_k\|}:1\leq k\leq n\}, \alpha_2=\min\{\frac{\|x_k\|^2-p_2^2+1}{2\|x_k\|}:1\leq k\leq n\},
\end{eqnarray*}
Then we have the inequality
\begin{eqnarray*}
(\alpha_1^2+\alpha_2^2)^\frac{1}{2}\sum_{k=1}^n \| x_k\|\leq\|\sum_{k=1}^n x_k\|
\end{eqnarray*}
where the equality holds if and only if
\begin{eqnarray*}
\sum_{k=1}^n x_k=(\alpha_1+i\alpha_2)\sum_{k=1}^n \|x_k\|a
\end{eqnarray*}

{\bf Proof.} From the first inequality in $(4)$ we have
$$\langle x_k-a,x_k-a\rangle\leq p_1^2$$
$$\| x_k\|^2+1-p_1^2\leq2Re\langle x_k, a\rangle, k=1,\cdots,n$$
$$\frac{\|x_k\|^2-p_1^2+1}{2\|x_k\|}\|x_k\|\leq Re\langle x_k,a\rangle$$
consequently
$$\alpha_1\|x_k\|\leq Re\langle x_k,a\rangle.$$
Similarly from the second inequality we obtain
$$\alpha_2\|x_k\|\leq Re\langle x_k,ia\rangle=Im\langle x_k,a\rangle.$$

Now apply Theorem 1 for $r_1=\alpha_1,r_2=\alpha_2.\Box$

{\bf Corollary 3.} Let $a$ be a unit vector in the complex inner product space $(H;\langle .,.\rangle)$. Suppose the vectors $x_k \in H-\{0\}, k\in\{1,\cdots,n\}$ such that 
$$ \|x_k-a\|\leq 1 , \|x_k-ia\|\leq 1.$$Then
\begin{eqnarray*}
\frac{\alpha}{\sqrt 2}\sum_{k=1}^n \| x_k\|\leq\|\sum_{k=1}^n x_k\|
\end{eqnarray*}
in which $ \alpha=\displaystyle{\min_{1\leq k\leq n}}\|x_k\|.$ The equality holds if and only if
\begin{eqnarray*}
\sum_{k=1}^n x_k=\alpha\frac{(1+i)}{2}\sum_{k=1}^n \|x_k\|a
\end{eqnarray*}

{\bf Proof.} Apply Theorem 2 for $\alpha_1=\frac{\alpha}{2}=\alpha_2.\Box$

{\bf Theorem 4.} Let $a$ be a unit vector in the inner product space $(H;\langle .,.\rangle)$ over the real or complex number field. Suppose that the vectors $x_k \in H-\{0\}, k\in\{1,\cdots,n\}$ satisfy 
$$ \|x_k-a\|\leq p , p\in(0,\sqrt{\alpha^2+1}),\alpha=\displaystyle{\min_{1\leq k\leq n}}\|x_k\|.$$

Then we have the inequality
$$\alpha_1\sum_{k=1}^n \|x_k\|\leq\|\sum_{k=1}^n x_k\|$$
where $$\alpha_1=\min\{\frac{\|x_k\|^2-p^2+1}{2\|x_k\|}:1\leq k\leq n\}$$.

The equality holds if and only if
$$\sum_{k=1}^n x_k=\alpha_1\sum_{k=1}^n \|x_k\|a$$

{\bf Proof.} The proof is similar to Theorem 2 in which we use Theorem 1 with $r_2=0.\Box$

The next Theorem is a generalization of Theorem 1. It is a modification of Theorem 3 of $\cite{DRA2}$, however our proof is apparently different.

{\bf Theorem 5.} Let $a_{1},\ldots ,a_{m}$ be orthonormal vectors in the complex inner product space $(H;\langle .,.\rangle)$. Suppose that for $1\leq t\leq m ,r_t,\rho_t\in R$ and that the vectors $x_k \in H, k\in\{1,\cdots,n\}$ satisfy \\
\begin{equation}
0\leq r_t ^2 \|x_k\|\leq Re\langle x_k,r_t a_t\rangle, 0\leq \rho_t ^2 \|x_k\| \leq Im\langle x_k,\rho_t a_t \rangle, t\in\{1,\cdots,m\}
\end{equation}
Then we have the inequality
\begin{equation}
(\sum_{t=1}^m r_t ^2+\rho_t ^2)^\frac{1}{2}\sum_{k=1}^n\|x_k\|\leq\|\sum_{k=1}^n x_k\|
\end{equation}
The equality holds in (7) if and only if
\begin{equation}
\sum_{k=1}^n x_k =\sum_{k=1}^n \|x_k\|\sum_{t=1}^m (r_t+i\rho_t)a_t.
\end{equation}

{\bf Proof.} If $\displaystyle{\sum_{t=1}^m}(r_t ^2 +\rho_t ^2)=0$, theorem is trivial. Assume that $\displaystyle{\sum_{t=1}^m}(r_t ^2 +\rho_t ^2)\neq 0.$ Summing inequalities (6) over $k$ from $1$ to $n$ and again over $t$ from $1$ to $m$ we get
$$\sum_{t=1}^m(r_t ^2 +\rho_t ^2)\sum_{k=1}^n\|x_k\|\leq Re\langle \sum_{k=1}^n x_k ,\sum_{t=1}^m r_t a_t \rangle+Im \langle \sum_{k=1}^n x_k ,\sum_{t=1}^m \rho_t a_t \rangle$$
$$=Re\langle \sum_{k=1}^n x_k ,\sum_{t=1}^m r_t a_t \rangle+Re \langle \sum_{k=1}^n x_k ,i\sum_{t=1}^m \rho_t a_t \rangle$$
$$= Re\langle \sum_{k=1}^n x_k ,\sum_{t=1}^m (r_t+i\rho_t)a_t \rangle$$
$$\leq |\langle \sum_{k=1}^n x_k ,\sum_{t=1}^m (r_t+i\rho_t)a_t \rangle|$$
$$\leq\|\sum_{k=1}^n x_k\|\|\sum_{t=1}^m (r_t+i\rho_t)a_t\|$$
$$=\|\sum_{k=1}^n x_k\|(\sum_{t=1}^m(r_t ^2 +\rho_t ^2))^\frac{1}{2}.$$
Then
$$(\sum_{t=1}^m(r_t ^2 +\rho_t ^2))^\frac{1}{2}\sum_{k=1}^n\|x_k\|\leq\|\sum_{k=1}^n x_k\|.$$
If $(8)$ holds, then
$$\|\sum_{k=1}^n x_k\|=\|\sum_{k=1}^n \|x_k\|\|\sum_{t=1}^m (r_t+i\rho_t)a_t\|=
\sum_{k=1}^n \|x_k\|(\sum_{t=1}^m(r_t ^2 +\rho_t ^2))^{1/2}.$$
Conversely, if the equality holds in $(7)$ we obtain from (6) that
$$(\sum_{t=1}^m(r_t ^2 +\rho_t ^2))^{1/2}\|\sum_{k=1}^n x_k\|=$$
$$=\sum_{t=1}^m(r_t ^2 +\rho_t ^2)\sum_{k=1}^n\|x_k\|\leq Re\langle \sum_{k=1}^n x_k ,\sum_{t=1}^m (r_t+i\rho_t)a_t \rangle\leq$$
$$\leq |\langle \sum_{k=1}^n x_k ,\sum_{t=1}^m (r_t+i\rho_t)a_t \rangle|\leq$$
$$\leq\| \sum_{k=1}^n x_k \|\|\sum_{t=1}^m (r_t+i\rho_t)a_t\|=\| \sum_{k=1}^n x_k \|(\sum_{t=1}^m(r_t ^2 +\rho_t ^2))^{1/2}$$ 
Thus we have
$$|\langle \sum_{k=1}^n x_k ,\sum_{t=1}^m (r_t+i\rho_t)a_t \rangle|=\| \sum_{k=1}^n x_k \|\|\sum_{t=1}^m (r_t+i\rho_t)a_t\|.$$
Consequently there exists $\eta\geq 0$ such that
$$ \sum_{k=1}^n x_k =\eta\sum_{t=1}^m (r_t+i\rho_t)a_t$$ from which we have
$$\eta(\sum_{t=1}^m (r_t^2+\rho_t^2))^\frac{1}{2}=\|\eta\sum_{t=1}^m (r_t+i\rho_t)a_t\|=\|\sum_{k=1}^n x_k\|=\sum_{k=1}^n\| x_k\|(\sum_{t=1}^m (r_t^2+\rho_t^2))^\frac{1}{2}$$Hence
$$\eta=\sum_{k=1}^n\| x_k\|.\Box$$

{\bf Corollary 6.} Let $a_{1},\ldots ,a_{m}$ be orthornormal vectors in the inner product space $(H;\langle .,.\rangle)$ over the real or complex number field. Suppose for $1\leq t\leq m $ that the vectors $x_k \in H, k\in\{1,\cdots,n\}$ satisfy 
$$0\leq r_t ^2 \|x_k\|\leq Re\langle x_k,r_t a_t\rangle.$$
Then we have the inequality
$$(\sum_{t=1}^m r_t ^2 )^\frac{1}{2}\sum_{k=1}^n\|x_k\|\leq\|\sum_{k=1}^n x_k\|.$$
The equality holds if and only if
$$\sum_{k=1}^n x_k =\sum_{k=1}^n \|x_k\|\sum_{t=1}^m r_t a_t.$$

{\bf Proof.} Apply Theorem 5 for $\rho_t=0.\Box$

{\bf Theorem 7.} Let $a_{1},\ldots ,a_{m}$ be orthornormal vectors in the complex inner product space $(H;\langle .,.\rangle)$. Suppose that the vectors $x_k \in H-\{0\}, k\in\{1,\cdots,n\}$ satisfy 
$$\|x_k-a_t\|\leq p_t ,\|x_k-ia_t\|\leq q_t, p_t, q_t\in(0,\sqrt{\alpha^2+1}),1\leq t\leq m$$ where $ \alpha=\displaystyle{\min_{1\leq k\leq n}}\|x_k\|.$ Let
$$ \alpha_t=\min\{\frac{\|x_k\|^2-p_t^2+1}{2\|x_k\|}:1\leq k\leq n\},\beta_t=\min\{\frac{\|x_k\|^2-q_t^2+1}{2\|x_k\|}:1\leq k\leq n\}.$$
Then we have the inequality
$$(\sum_{t=1}^m\alpha_t^2+\beta_t^2)^\frac{1}{2}\sum_{k=1}^n\|x_k\|\leq\|\sum_{k=1}^n x_k\|$$
where equality holds if and only if
$$\sum_{k=1}^n x_k=\sum_{k=1}^n\|x_k\|\sum_{t=1}^m(\alpha_t+i\beta_t)a_t.$$

{\bf Proof.} For $1\leq t\leq m, 1\leq k\leq n$ it follows from $\|x_k-a_t\|\leq p_t$ that
$$\langle x_k-a_t\rangle,x_k-a_t\rangle\leq p_t^2$$
$$\frac{\|x_k\|^2-p_t^2+1}{2\|x_k\|}\|x_k\|\leq Re\langle x_k,a_t\rangle0$$
$$\alpha_t\|x_k\|\leq Re\langle x_k,a_t\rangle$$
and similarly
$$\beta_t\|x_k\|\leq Re\langle x_k,ia_t\rangle=Im\langle x_k,a_t\rangle, $$
Now applying Theorem 4 with $r_t=\alpha_t, \rho_t=\beta_t$ we deduce the desired inequality.$\Box$

{\bf Corollary 8.} Let $a_{1},\ldots ,a_{m}$ be orthornormal vectors in the complex inner product space $(H;\langle .,.\rangle)$. Suppose that the vectors $x_k \in H, k\in\{1,\cdots,n\}$ satisfy 
$$ \|x_k-a_t\|\leq 1 ,\|x_k-ia_t\|\leq 1, 1\leq t\leq m $$ 
Then 
$$\frac{\alpha}{\sqrt{2}}\sqrt{m}\sum_{k=1}^n\|x_k\|\leq\|\sum_{k=1}^n x_k\|.$$ 
The equality holds if and only if $$\sum_{k=1}^n x_k=\alpha\frac{(1+i)}{2}\sum_{k=1}^n\|x_k\|\sum_{t=1}^m a_t$$

{\bf Proof.} Applying Theorem 7 for $\alpha_t=\frac{\alpha}{2}=\beta_t.\Box$

{\bf Remark.} It is interesting to note that

$$\frac{\alpha}{\sqrt{2}}\sqrt{m}\leq\frac{\|\sum_{k=1}^n x_k\|}{\sum_{k=1}^n\|x_k\|}\leq 1.$$
$$\alpha\leq\sqrt{\frac{2}{m}}$$

{\bf Corollary 9.} Let $a$ be a unit vector in the complex inner product space $(H;\langle .,.\rangle)$. Suppose that the vectors $x_k \in H-\{0\}, k\in\{1,\cdots,n\}$ satisfy \\
$$ \|x_k-a\|\leq p_1, \|x_k-ia\|\leq p_2, p_1, p_2\in(0,1].$$ Let 
$$ \alpha_1=\min\{\frac{\|x_k\|^2-p_1^2+1}{2\|x_k\|}:1\leq k\leq n\},\alpha_2=\min\{\frac{\|x_k\|^2-p_2^2+1}{2\|x_k\|}:1\leq k\leq n\}.$$ 
If $\alpha_1\neq(1-p_1^2)^\frac{1}{2}$, or $\alpha_2\neq(1-p_2^2)^\frac{1}{2}$, then we have the following strictly inequality 
$$(2-p_1^2-p_2^2)^\frac{1}{2}\sum_{k=1}^n \| x_k\|<\|\sum_{k=1}^n x_k\|$$

{\bf Proof.} If equality holds, then by Theorem 2 we have 
$$(\alpha_1^2+\alpha_2^2)^\frac{1}{2}\sum_{k=1}^n \| x_k\|\leq\|\sum_{k=1}^n x_k\|=(2-p_1^2-p_2^2)^\frac{1}{2}\sum_{k=1}^n \| x_k\|$$ and so
$$(\alpha_1^2+\alpha_2^2)^\frac{1}{2}\leq(2-p_1^2-p_2^2)^\frac{1}{2}.$$
On the other hand for $1\leq k\leq n$, 
$$\frac{\|x_k\|^2-p_1^2+1}{2\|x_k\|}\geq(1-p_1^2)^\frac{1}{2}$$ and so
$$\alpha_1 \geq (1-p_1^2)^\frac{1}{2}.$$
Similarly
$$\alpha_2\geq(1-p_2^2)^\frac{1}{2}.$$
Hence
$$(2-p_1^2-p_2^2)^\frac{1}{2}\leq(\alpha_1^2+\alpha_2^2)^\frac{1}{2}$$ Thus
$$\sqrt{\alpha_1^2+\alpha_2^2}=(2-p_1^2-p_2^2)^\frac{1}{2}.$$
Therefore
$$\alpha_1=(1-p_1^2)^\frac{1}{2} {\rm and} \alpha_2=(1-p_2^2)^\frac{1}{2}$$
a contradiction.$\Box$

The following result looks like Corollary 2 of $\cite{DRA2}$.

{\bf Theorem 10}. Let $a$ be a unit vector in the complex inner product space $(H;\langle .,.\rangle), M\geq m >0, L\geq \ell>0$ and $x_k \in H-\{0\}, k\in\{1,\cdots,n\}$ such that
$$Re\langle Ma-x_k,x_k-ma\rangle\geq0,Re\langle L ia-x_k,x_k-\ell ia\rangle\geq0$$ or equivalently,
$$\|x_k-\frac{m+M}{2}a\|\leq\frac{M-m}{2},\|x_k-\frac{L+\ell}{2}ia\|\leq\frac{L-\ell}{2}.$$
Let
$$\alpha_{m,M}=\min\{\frac{\|x_k\|^2+mM}{(m+M)\|x_k\|}:1\leq k\leq n\}$$ and $$\alpha_{\ell,L}=\min\{\frac{\|x_k\|^2+\ell L}{(\ell+L)\|x_k\|}:1\leq k\leq n\}$$ Then we have the inequlity
$$(\alpha_{m,M}^2+\alpha_{\ell,L}^2)^\frac{1}{2}\sum_{k=1}^n\|x_k\|\leq\|\sum_{k=1}^n x_k\|.$$
The equality holds if and only if
$$\sum_{k=1}^n x_k=(\alpha_{m,M}+i\alpha_{\ell,L})\sum_{k=1}^n\|x_k\|a.$$

{\bf Proof}. For each $1\leq k\leq n$, it follows from 
$$\|x_k-\frac{m+M}{2}a\|\leq\frac{M-m}{2}$$ 
that
$$\langle x_k-\frac{m+M}{2}a, x_k-\frac{m+M}{2}\rangle\leq(\frac{M-m}{2})^2.$$ 
Hence
$$\|x_k\|^2+mM\leq(m+M)Re\langle x_k, a\rangle. $$
Then
$$\frac{\|x_k\|^2+mM}{(m+M)\|x_k\|}\|x_k\|\leq Re\langle x_k, a\rangle$$
consequently
$$\alpha_{m,M}\|x_k\|\leq Re\langle x_k, a\rangle.$$
Similarly from the second inequlity we deduce
$$\alpha_{\ell,L}\|x_k\|\leq Im\langle x_k, a\rangle.$$ Applying Theorem 1 for $ r_1=\alpha_{m,M},r_2=\alpha_{\ell,L}$, we infer the desired inequality.$\Box$

{\bf Theorem 11}. Let $a$ be a unit vector in the complex inner product space $(H;\langle .,.\rangle), M \geq m > 0, L\geq \ell > 0$ and $x_k \in H-\{0\}, k\in\{1,\cdots,n\}$ such that 
$$Re\langle Ma-x_k,x_k-ma\rangle\geq 0,Re\langle Lia-x_k,x_k-\ell ia\rangle\geq 0$$
or equivalently
$$\|x_k-\frac{m+M}{2}a\|\leq\frac{M-m}{2},\|x_k-\frac{L+\ell}{2}ia\|\leq\frac{L-\ell}{2}.$$Let
$$\alpha_{m,M}=\min\{\frac{\|x_k\|^2+mM}{(m+M)\|x_k\|}:1\leq k\leq n\}$$
and
$$\alpha_{\ell,L}=\min\{\frac{\|x_k\|^2+\ell L}{(\ell+L)\|x_k\|}:1\leq k\leq n\}.$$
If $\alpha_{m,M}\neq 2\frac{\sqrt{mM}}{m+M}$, or $\alpha_{\ell,L}\neq 2\frac{\sqrt{\ell L}}{\ell+L}$, then we have
$$2(\frac{mM}{(m+M)^2}+\frac{\ell L}{(\ell+L)^2})^\frac{1}{2}\sum_{k=1}^n\|x_k\|<\|\sum_{k=1}^n x_k\|.$$

{\bf Proof}. If $2(\frac{mM}{(m+M)^2}+\frac{\ell L}{(\ell+L)^2})^\frac{1}{2}\displaystyle{\sum_{k=1}^n}\|x_k\|=\|\displaystyle{\sum_{k=1}^n} x_k\|$ then by theorem 10 we have
$$(\alpha_{m,M}^2+\alpha_{\ell,L}^2)^\frac{1}{2}\sum_{k=1}^n\|x_k\|\leq\|\sum_{k=1}^n x_k\|\\
=2(\frac{mM}{(m+M)^2}+\frac{\ell L}{(\ell+L)^2})^\frac{1}{2}\sum_{k=1}^n\|x_k\|.$$
Consequently
$$(\alpha_{m,M}^2+\alpha_{\ell,L}^2)^\frac{1}{2}\leq 2(\frac{mM}{(m+M)^2}+\frac{\ell L}{(\ell+L)^2})^\frac{1}{2}.$$
On the other hand for $ 1\leq k\leq n,\frac{\|x_k\|^2+mM}{(m+M)\|x_k\|}\geq 2\frac{\sqrt{mM}}{m+M}$, and $\frac{\|x_k\|^2+\ell L}{(\ell+L)\|x_k\|}\geq 2\frac{\sqrt{\ell L}}{\ell+L}$, so
$$(\alpha_{m,M}^2+\alpha_{\ell,L}^2)^\frac{1}{2}\geq 2(\frac{mM}{(m+M)^2}+\frac{\ell L}{(\ell+L)^2})^\frac{1}{2}.$$
Then
$$(\alpha_{m,M}^2+\alpha_{\ell,L}^2)^\frac{1}{2}=2(\frac{mM}{(m+M)^2}+\frac{\ell L}{(\ell+L)^2})^\frac{1}{2}.$$
Hence
$$\alpha_{m,M}= 2\frac{\sqrt{mM}}{m+M}$$ 
and
$$ \alpha_{\ell,L}= 2\frac{\sqrt{\ell L}}{\ell+L}$$
a contradection.$\Box$

Finally we mention two applications of our results to the complex numbers.

{\bf Corollary 12.} Let $a\in C$ with $|a|=1$. Suppose that $z_k \in C, k\in\{1,\cdots,n\}$ such that 
$$ |z_k-a|\leq p_1 , |z_k-ia|\leq p_2,p_1, p_2\in(0,\sqrt{\alpha^2+1})$$ 
where $$\alpha=\min\{|z_k|:1\leq k\leq n\}.$$
Let 
$$ \alpha_1=\min\{\frac{|z_k|^2-p_1^2+1}{2|z_k|}:1\leq k\leq n\},\alpha_2=\min\{\frac{|z_k|^2-p_2^2+1}{2|z_k|}:1\leq k\leq n\}.$$
Then we have the inequality
$$\sqrt{\alpha_1^2+\alpha_2^2}\sum_{k=1}^n | z_k|\leq|\sum_{k=1}^n z_k|.$$
The equality holds if and only if
$$\sum_{k=1}^n z_k=(\alpha_1+i\alpha_2)(\sum_{k=1}^n |z_k|)a.$$

{\bf Proof.} Apply Theorem 2 for $H=C.\Box$

{\bf Corollary 13.} Let $a\in C$ with $|a|=1$. Suppose that $z_k \in C, k\in\{1,\cdots,n\}$ such that 
$$ |z_k-a|\leq 1 , |z_k-ia|\leq 1. $$ 
If $\alpha=\min\{|z_k|:1\leq k\leq n\}$. Then we have the inequality
$$\frac{\alpha}{\sqrt 2}\sum_{k=1}^n |z_k|\leq|\sum_{k=1}^n z_k|$$
the equality holds if and only if
$$\sum_{k=1}^n z_k=\alpha\frac{(1+i)}{2}(\sum_{k=1}^n |z_k|)a.$$

{\bf Proof.} Apply Corollary 3 for $H=C.\Box$

\end{document}